\documentclass[11pt]{article}

\usepackage[margin=1in]{geometry}
\usepackage{fix-cm}
\usepackage{graphicx}
\usepackage{amsmath,amsfonts,amssymb}
\usepackage{newtxtext}
\usepackage{newtxmath}
\usepackage[round]{natbib}
\usepackage{hyperref}
\usepackage{placeins}
\usepackage{bm}
\usepackage{tikz}
\usepackage{url}

\graphicspath{{figs/}}

\DeclareMathOperator{\tr}{tr}
\newcommand{\dd}{\partial}
\newcommand{\norm}[1]{\lVert #1\rVert}
\newcommand{\aT}{\boldsymbol{a}}

\newcommand{\dt}{\delta^{+}}

\title{Near-cancellation of enstrophy and strain requires shear or a balanced vortex}

\author{Jonathan Massey\\
Center for Turbulence Research, Stanford University, Stanford, CA 94305, USA\\
\texttt{masseyj@stanford.edu}}

\date{}

\begin{document}

\maketitle

\begin{abstract}
In the nonlinear source of the pressure Poisson equation the fluctuating enstrophy and strain nearly cancel. We ask where the cancellation is complete, which we call silent, and how the source escapes it. Pointwise the source depends only on the eigenvalues of the velocity gradient. Pure shear has none, so it is silent. It is the origin of the normalised invariant plane, whose second-invariant coordinate is the normalised source. Let $f$ be the shear fraction of the gradient and $m$ its distance from that origin. We prove that $m\le1-f\le3m$ at every point of every incompressible flow. Wherever the source nearly vanishes and the third invariant is small, the gradient is nearly pure shear. Otherwise the cancelling gradient is a balanced vortex, a swirl its own strain cancels at any strength. At a no-slip wall the silent object is exact: the wall shear-stress fluctuation, a rank-one sheet. Away from the wall we assume the same sheet, dressed with a small residual, and derive how it escapes. The escape passes through one element of the residual, the along-sheet variation of the sheet-normal velocity. Its amplitude is the geometric mean of sheet and residual, and its sign selects swirl or strain. A rank-two silent object has the same source through the same gates but balanced eigenvalues, so the assumption fixes the escape's topology, not its source. The shear fraction is thus a pointwise proxy for cancellation by shear, and a model of the source must be built on what escapes.
\end{abstract}

\section{Introduction}\label{sec:intro}

At a wall, the fluctuating kinematic pressure, hereafter `pressure', is the load and sources the sound through scattering and radiation. Within the flow, its correlation with the strain rate is the redistribution term that second-moment closures model. Physically, the pressure is the Lagrange multiplier that enforces incompressibility. Its gradient projects the momentum equation onto the manifold of divergence-free velocity fields, and the source of the resulting Poisson equation is the divergence the advection term would otherwise create. In a large-eddy simulation the unresolved part of the source, the double divergence of the subgrid stress, is what the subgrid model must supply.

For incompressible flow the divergence of the momentum equation leaves a Poisson equation for the pressure. \citet{bradshaw_note_1981} pointed out that its right-hand side is the difference of two squares. For the instantaneous velocity $\tilde u_i$, with vorticity $\tilde\omega_i$ and strain rate $\tilde s_{ij}$, $\dd_i\dd_i\tilde p=\tfrac12\tilde\omega_i\tilde\omega_i-\tilde s_{ij}\tilde s_{ij}$. Under Reynolds decomposition the fluctuating pressure is forced by two sources \citep{mansour_reynolds-stress_1988,kim_structure_1989}.

The linear source is proportional to the mean shear. The nonlinear source is the turbulence acting on itself. Writing $a_{ij}=\dd_ju_i$ for the fluctuating gradient and $\aT$ for the tensor,
\begin{equation}
    \dd_i\dd_i p=-2\,\dd_j\overline{u}_i\,\dd_i u_j-\bigl[\tr(\aT^2)-\langle\tr(\aT^2)\rangle\bigr],
    \qquad -\tr(\aT^2)=\tfrac12\omega^2-s^2,
    \label{eq:split}
\end{equation}
with $\omega^2=\omega_i\omega_i$ and $s^2=s_{ij}s_{ij}$ the fluctuating enstrophy and strain, and $a\equiv\norm{\aT}=(\tfrac12\omega^2+s^2)^{1/2}$ the gradient magnitude. Where the gradient is intense and the source is small, the nonlinear source is the small remainder of two large terms. What such a gradient must look like is the question of \S\ref{sec:bounds}.

For channels, pipes and zero-pressure-gradient boundary layers the inner-scaled wall-pressure variance grows with the friction Reynolds number $\dt\equiv u_\tau\delta/\nu$ \citep{klewicki_statistical_2008,panton_correlation_2017}. Here $u_\tau$ is the friction velocity, $\nu$ the kinematic viscosity and $\delta$ the outer length: the channel half-height, the pipe radius or the boundary-layer thickness. Evidence suggests the growth belongs to the nonlinear source since the inner-scaled linear contribution saturates to invariance \citep{massey_linear_2026}. Its proposed carrier is the sheet-like vortical fissures that bound the uniform-momentum zones of the inertial layer \citep{meinhart_existence_1995,de_silva_interfaces_2017}. \citet{bautista_uniform_2019} propose a probabilistic skeletal model of the velocity field, built from the inertial-layer velocity jumps across these fissures. Fed through the Poisson equation, it locates the source on the fissures, as the product of the fissure's shear and the wall-normal velocity residual across it \citep{massey_behind_2026}. The fissure's own shear is the largest gradient in the flow, yet on its own it sources nothing and enters only multiplied by the residual. Why this is so is a question about the algebra of the source. That algebra holds at every point of every incompressible flow, and it is the subject of this paper.

Part of the answer is established. The pure-shear constituent of the velocity gradient cancels exactly in the second invariant \citep{keylock_schur_2018,das_revisiting_2020}. The converse is not. Bounds on the departure from normality, one-sided \citep{henrici_bounds_1962} and two-sided \citep{lee_best_1996}, are written in norms and commutators, not in the invariants. The normalised-gradient framework of \citet{das_reynolds_2019} bounds the realisable invariant plane, but not the shear content at each point of it. The Lee bound of \citet{keylock_lee_2019} is a one-sided upper limit on the non-normality a swirling region can attain. This paper proves the converse as the two-sided wedge $m\le1-f\le3m$ of \S\ref{sec:bounds}. The ingredients are classical: the Schur identity above, and a bound on the roots of a cubic by its coefficients. What is new is the statement in the invariants of turbulence, with its constants, and the one family that evades it. On the wedge we build a mechanism. Complete cancellation lives on sheets of viscous shear-stress fluctuation, exactly so at the wall, and \S\ref{sec:escape} derives the conditions under which a sheet escapes it. The square-root splitting of the eigenvalues of a perturbed nilpotent is also classical \citep{moro_lidskii_1997}. What \S\ref{sec:escape} adds is the identification of the splitting parameter with one element of the residual in the frame of the sheet, and the sign rule that follows. The mechanism is a hypothesis about turbulence. It rests on two assumptions, and \S\ref{sec:conclusion} says what would test them and what follows for the wall pressure and for models.

\section{Invariant framework}\label{sec:identity}

The real, trace-free fluctuating gradient tensor $\aT(\boldsymbol x,t)$ has eigenvalues $\lambda_1,\lambda_2,\lambda_3$, and its invariants are $Q_a=-\tfrac12\tr(\aT^2)$ and $R_a=-\det\aT$ \citep{chong_general_1990}. In Schur form, $\aT=\mathsf U(\Lambda+\mathsf N)\mathsf U^{*}$ with $\mathsf U$ unitary, $\Lambda$ diagonal and $\mathsf N$ strictly upper triangular. The norm of $\mathsf N$ is Henrici's departure from normality, which \citet{keylock_schur_2018} introduced to turbulence as the frame-invariant shearing content. Since $\Lambda\mathsf N$, $\mathsf N\Lambda$ and $\mathsf N^2$ are strictly triangular, $\tr(\aT^2)=\sum_i\lambda_i^2$, while unitary invariance gives $a^2=\sum_i|\lambda_i|^2+\norm{\mathsf N}^2$. Hence, with $\sigma^2\equiv\norm{\mathsf N}^2$,
\begin{equation}
    \tfrac12\omega^2-s^2=-\sum_i\lambda_i^2,\qquad
    \tfrac12\omega^2=\tfrac12\sigma^2+\sum_i(\mathrm{Im}\,\lambda_i)^2,\qquad
    s^2=\tfrac12\sigma^2+\sum_i(\mathrm{Re}\,\lambda_i)^2 .
    \label{eq:WSsplit}
\end{equation}
The shearing content enters enstrophy and strain in equal halves and drops out of the source. The source is the eigenvalue field. Swirl sources $\dd_i\dd_ip>0$ and normal straining sources $\dd_i\dd_ip<0$. The shearing content, at fixed eigenvalues, sources nothing. Hereafter we call a gradient, or the part of one, \emph{silent} when it contributes nothing to the source.\footnote{The triple decomposition of \citet{kolar_vortex_2007} is the standardised real Schur form \citep{kronborg_triple_2023}. \citet{das_revisiting_2020} write the pressure Laplacian in it, retaining a shear--rotation cross term, and find shear the largest constituent of $a^2$. It partitions each complex-pair block differently, so its residual does not carry the eigenvalues. The exact statement is \eqref{eq:WSsplit}.} We need four dimensionless measures, defined where $a^2>0$:
\begin{equation}
    \chi=\frac{|\tfrac12\omega^2-s^2|}{a^2},\qquad
    \hat q=\frac{Q_a}{a^2},\qquad
    \hat r=\frac{R_a}{a^3},\qquad
    f=\frac{\sigma^2}{a^2}.
    \label{eq:defs}
\end{equation}
These are the cancellation defect $\chi$, the coordinates $(\hat q,\hat r)$ of the normalised invariant plane \citep{das_reynolds_2019}, and the shear fraction $f$, which is unity for a simple shearing motion. That nilpotent gradient is what we call pure shear throughout.

\section{Cancellation and shear}\label{sec:bounds}

\subsection{The wedge}\label{sec:wedge}

Combining \eqref{eq:WSsplit} with the definitions \eqref{eq:defs}, shear implies cancellation, with constants, in three parts:
\begin{equation}
    \chi\le 1-f,\qquad
    \min(\tfrac12\omega^2,s^2)\ge\tfrac12 fa^2,\qquad
    |\hat r|\le\Bigl(\frac{1-f}{3}\Bigr)^{3/2}.
    \label{eq:forward}
\end{equation}
The first is $|\sum\lambda_i^2|\le\sum|\lambda_i|^2=(1-f)a^2$. The second is read off \eqref{eq:WSsplit}. The third is the arithmetic--geometric mean inequality applied to $|\!\det\aT|=\prod|\lambda_i|$. The converse is the statement that a point close to the origin of the invariant plane has little eigenvalue content. It is a Cauchy-type bound on the roots of the characteristic cubic $\lambda^3+Q_a\lambda+R_a=0$ by its coefficients \citep{wimmer_polynomials_2015}. Each eigenvalue satisfies $|\lambda|^3\le|Q_a||\lambda|+|R_a|$. If $|\lambda|$ exceeded both $(2|Q_a|)^{1/2}$ and $(2|R_a|)^{1/3}$, the right-hand side would be smaller than $|\lambda|^3$. So each $|\lambda_i|^2\le ma^2$ with
\begin{equation}
    m\equiv\max\bigl(2|\hat q|,\,(2|\hat r|)^{2/3}\bigr),
    \qquad\text{and summing,}\qquad
    1-f\le 3m .
    \label{eq:converse}
\end{equation}
The first and third parts of \eqref{eq:forward} give $1-f\ge m$ in either branch of the maximum. The two bounds therefore close into a wedge,
\begin{equation}
    m\le 1-f\le 3m .
    \label{eq:wedge}
\end{equation}
This is the central algebraic result of the paper. Up to a factor of three, $1-f$ is $m$, a distance from the pure-shear origin of the invariant plane in the metric of \eqref{eq:converse}. The lower edge is attained by every tensor whose spectrum is real, or purely imaginary. The upper edge is not attained. Both $1-f$ and $m$ are $a^{-2}$ times a function of the spectrum, so their ratio depends only on the spectrum's shape. On the compact set of normalised spectra with $m>0$ the ratio is continuous and takes its maximum. Write the spectrum as $\mu(1\pm\mathrm i\sqrt u),\,-2\mu$. On the $\hat r$ branch the ratio falls and then rises with $u$, and on the $\hat q$ branches it is monotone towards the crossings of the two branches, so the maximum lies at a crossing. The branches cross twice. At $u=1$ the ratio is exactly $2$. At the other crossing, on the swirl side, $u$ is the real root of $(u-3)^3=2(1+u)^2$, about $8.75$, and the ratio is $(u+3)/(u-3)$, which is $2.04$ to three figures. Appendix~\ref{app:tests} shows the wedge and this constant on a random-tensor ensemble.

Thus, wherever enstrophy and strain cancel to a tolerance $\eta$, that is $\chi\le\eta$ and $|\hat r|\le\tfrac12\eta^{3/2}$, the gradient is shear-dominated, $f\ge1-3\eta$. The wedge is also a realisability constraint on $(\hat q,\hat r,f)$, loosely as the anisotropy triangle constrains the Reynolds stresses.

\subsection{The shear-dominated gradient}\label{sec:shear}

The shear fraction is unity if and only if all three eigenvalues vanish, that is, if $\aT$ is nilpotent. The distance from $\aT$ to the nearest real trace-free nilpotent is at most $\sqrt{1-f}\,a$. In the standardised real Schur form \citep{kronborg_triple_2023} each complex pair occupies a block with equal diagonal entries. Zeroing the diagonal, and the smaller off-diagonal entry of each such block, removes a norm of at most $\sum_i|\lambda_i|^2$, and less when a pair is complex, and leaves a real nilpotent. A non-zero nilpotent is rank one or rank two. The rank-one nilpotent is $\aT=g\,\boldsymbol t\otimes\boldsymbol n$ with $\boldsymbol t\perp\boldsymbol n$, which is locally a plane shear layer. Here $\boldsymbol n$ is the unit normal to the layer, $\boldsymbol t$ is the unit vector along the velocity that varies across it, and $g$ is the shear, the velocity gradient across the layer. On the pure layer $g=a$. The rank-two nilpotent is the superposition of two such layers, and rank two is the generic case. Which of the two a gradient is near is measured by $r_1=1-s_1^2/a^2$, with $s_1$ the largest singular value. It vanishes on a rank-one nilpotent and is one half on the equal-strength rank-two one. In the rank-one case the strain eigenvalues are $(g/2,0,-g/2)$, and the vorticity has magnitude $g$ and lies along the intermediate principal axis. The fluctuating viscous stress transmitted across the layer is purely tangential, $2\nu s_{ij}n_j=\nu g\,t_i$. A rank-one silent gradient is therefore a fluctuating viscous shear stress with matched rotation. It is the local form of the small-scale shear layers of \citet{watanabe_response_2023}. Its exact alignment of vorticity with the intermediate strain axis is consistent with the statistical preference \citet{ashurst_alignment_1987} found.

\subsection{The exception}\label{sec:exception}

Exact cancellation, $\chi=0$, occurs if and only if the spectrum is equilateral, $\{\lambda_i\}=\lambda_0\{1,\mathrm e^{\pm2\pi\mathrm i/3}\}$. For $\lambda_0=0$ this is the nilpotent, or shear, case. For $\lambda_0\neq0$ it is a one-parameter family along $\hat q=0$: any amount of shear on a critically strained vortex. The shear-free member sits at $|\hat r|=3^{-3/2}$, the corner of the realisable invariant plane of \citet{das_reynolds_2019}, which the third part of \eqref{eq:forward} re-derives in one line. Along the family $1-f=3|\hat r|^{2/3}$, so the third-invariant condition of \S\ref{sec:wedge} is necessary. There are two ways of cancelling: a shearing motion with no eigenvalue content, and a balanced vortex. Which of them turbulence uses, and where, is a question for data.

\section{Escape from cancellation}\label{sec:escape}

\subsection{The wall}\label{sec:wall}

At a smooth no-slip wall the silence is exact. No slip removes every wall-parallel derivative at $y=0$, and continuity then removes $\dd_2u_2$. With $\boldsymbol\tau=\nu(\dd_2u_1,0,\dd_2u_3)|_{y=0}$ the fluctuating wall shear-stress vector, kinematic like the pressure,
\begin{equation}
    \aT\big|_{y=0}=\frac1\nu\,\boldsymbol\tau\otimes\boldsymbol e_2,
    \qquad\text{hence}\qquad
    f=1,\quad\chi=0,\quad \tfrac12\omega^2=s^2=\frac{|\boldsymbol\tau|^2}{2\nu^2}\ \ \text{identically}.
    \label{eq:wall}
\end{equation}
The wall plane is an exact-cancellation sheet, and its intensity field is the wall shear-stress fluctuation itself. It is also escape-proof to leading order. Expanding about $y=0$, no slip and continuity make $u_2=O(y^2)$, so the wall-parallel derivatives of $u_2$ carry no $O(y)$ term. Hence $\tr(\aT^2)=O(y^2)$ while $a^2=|\boldsymbol\tau|^2/\nu^2+O(y)$. The expansion holds in the viscous sublayer. That $\tfrac12\omega^2=s^2$ at a wall is classical, and so is the form \eqref{eq:wall}. \citet{chen_boundary_2024} derive the nilpotent boundary value in general, with the vanishing of the boundary flux of the third invariant. What we take from it is the rank-one nilpotency imposed kinematically by no slip, whose near-wall budget of normal and non-normal parts \citet{keylock_role_2025} supplies.

\subsection{Away from the wall}\label{sec:sheet}

\subsubsection{The assumptions}

Away from the wall no kinematic constraint fixes the silent object. Its form is a conditional tendency of the flow, not an identity, and we must assume it. We make two assumptions. The first is that the silent object is rank one, as the wall forces and as the internal shear layers and fissures of the literature suggest, rather than the generic rank-two nilpotent of \S\ref{sec:shear}. The second is that everything else in the gradient is small beside it. The object is then a thin sheet dressed with a residual, $\aT=g\,\boldsymbol t\otimes\boldsymbol n+\epsilon\boldsymbol b$. The sheet is the rank-one layer of \S\ref{sec:shear}, with normal $\boldsymbol n$, velocity direction $\boldsymbol t$ and shear $g$. The residual $\epsilon\boldsymbol b$ has $\boldsymbol b$ a trace-free tensor of unit norm and $\epsilon$ its amplitude, with $g\gg\epsilon$. Section~\ref{sec:conclusion} suggests how both assumptions could be tested. In the inertial layer it is a vortical fissure \citep{bautista_uniform_2019,massey_behind_2026}. Since $\tr[(\boldsymbol t\otimes\boldsymbol n)^2]=0$,
\begin{equation}
    \tr(\aT^2)=2g\epsilon\,n_ib_{ij}t_j+\epsilon^2\tr(\boldsymbol b^2).
    \label{eq:fissure}
\end{equation}
The invariants are $Q_a=-g\epsilon\,n_ib_{ij}t_j+O(\epsilon^2)$ and $R_a=O(g\epsilon^2)$, with no term linear in $\epsilon$ in the determinant since a rank-one tensor has no cofactors, and the characteristic equation $\lambda^3+Q_a\lambda+R_a=0$ has two escaping roots and one small one:
\begin{equation}
    \lambda\simeq\pm\bigl(g\epsilon\,n_ib_{ij}t_j\bigr)^{1/2}.
    \label{eq:escape}
\end{equation}
The pair is real when the coupling $n_ib_{ij}t_j$ is positive and imaginary when it is negative. The third root is $-R_a/Q_a$ to leading order, of order $\epsilon$, smaller than the pair by $(\epsilon/g)^{1/2}$ at order-one coupling. The square-root splitting itself is the classical response of a defective eigenvalue to perturbation \citep{moro_lidskii_1997}. Equation~\eqref{eq:escape} names the parameter that splits it, and it separates three things that are usually run together.

\subsubsection{The gate}

A sheet escapes silence through the single matrix element $n_ib_{ij}t_j$ of its residual. With $\boldsymbol u_b$ the residual velocity, this element is $\dd_t(\boldsymbol u_b\cdot\boldsymbol n)$, the along-sheet variation of the velocity the residual carries across the sheet. In wall coordinates, with $\boldsymbol n=\boldsymbol e_2$, a sheet with streamwise $\boldsymbol t=\boldsymbol e_1$ has the gate $\dd_1u_2$ and the source $-2\,\dd_2u_1\,\dd_1u_2$. A sheet with spanwise $\boldsymbol t=\boldsymbol e_3$ has the gate $\dd_3u_2$ and the source $-2\,\dd_2u_3\,\dd_3u_2$. \citet{kim_structure_1989} and \citet{chang_relationship_1999} found the second to be the leading cross term of the nonlinear source in the buffer layer of low-Reynolds-number channel flow. In this frame it is the gate of a sheet whose velocity direction is spanwise. Equation~\eqref{eq:escape} is the cross term written in the frame of the sheet, whatever its orientation. Without it the sheet contributes nothing at any strength, and only the residual's own source, of order $\epsilon^2$, remains. Escape is a condition on alignment, not on amplitude.

\subsubsection{The amplification}

What does escape is raised to the geometric mean $(g\epsilon)^{1/2}$ of the sheet's silent strength and the residual's own. The source is the square of that, $-\tr(\aT^2)=-2g\epsilon\,n_ib_{ij}t_j$ to leading order, linear in the sheet strength. The shear-stress fluctuation sets the amplitude of everything that escapes while contributing none of it.

\subsubsection{The sign}

The sign of the coupling decides what the escape is. In the plane of $\boldsymbol t$ and $\boldsymbol n$ the sheet's rotation is $\tfrac12(g-\epsilon c)$ and its shear strain $\tfrac12(g+\epsilon c)$, with $c=n_ib_{ij}t_j$ and the residual's own $t_ib_{ij}n_j$ absorbed into $g$. A negative coupling raises the rotation over the strain and gives an imaginary pair: swirl, sourcing $\dd_i\dd_ip>0$, the local streamlines closing. It does so at any strength that leaves the ridge, since a real trace-free tensor with $Q_a>0$ always has a complex pair. A positive coupling gives a real pair: strain, sourcing $\dd_i\dd_ip<0$, the local streamlines hyperbolic. It does so only where the coupling is strong enough for the second invariant to dominate the cubic, $|c|\gg(\epsilon/g)^{1/3}$. Below that the source is negative but the pair is still complex, a strain-dominated focus. For a sheet $u_1=gx_2$ dressed with $u_2=\epsilon\sin kx_1$ the gate alternates along the sheet, and swirl and strain alternate as the cores and braids of a rolling-up sheet. This is the frozen local topology, not a statement about material motion. Nothing in the algebra prefers one sign of the source. The topology is asymmetric: swirl is unconditional and strain is not.

\subsubsection{What the rank assumption fixes}

\begin{figure}
    \centering
    \begin{tikzpicture}[>=stealth,scale=1.15]
        \draw[->] (-3.4,0) -- (3.4,0) node[right] {$\hat r$};
        \draw[->] (0,-2.1) -- (0,2.1) node[above] {$\hat q$};
        \draw[very thick] (-2.6,0) -- (2.6,0);
        \fill (-2.6,0) circle (1.6pt);
        \fill (2.6,0) circle (1.6pt);
        \node[above] at (2.6,0.08) {$3^{-3/2}$};
        \node[above] at (-2.6,0.08) {$-3^{-3/2}$};
        \fill (0,0) circle (2.2pt);
        \node[above left] at (-0.08,0.12) {pure shear, $f=1$};
        \draw[->,thick] (0,0.15) -- (0,1.6) node[midway,right] {sheet $+$ residual: swirl, $\dd_i\dd_ip>0$};
        \draw[->,thick] (0,-0.15) -- (0,-1.6) node[midway,right] {sheet $+$ residual: strain, $\dd_i\dd_ip<0$};
        \draw[->,thick,dashed] (-0.12,-0.12) -- (-1.3,-1.3) node[below left] {rank two $+$ residual};
        \node[below] at (-1.9,-0.1) {balanced vortices, silent};
    \end{tikzpicture}
    \caption{The normalised invariant plane near its origin. Pure shear is the origin and is silent. A dressed rank-one sheet leaves it along $\hat q$ and sources pressure, swirl upward and strain downward. A dressed rank-two nilpotent leaves it with $\hat q$ and $\hat r$ of the same order, in any quadrant, with balanced eigenvalues and a source of the same order as the sheet's. The ridge $\hat q=0$ is the balanced vortices, silent at any strength, ending at the corners $|\hat r|=3^{-3/2}$ of the realisable plane.}
    \label{fig:plane}
\end{figure}
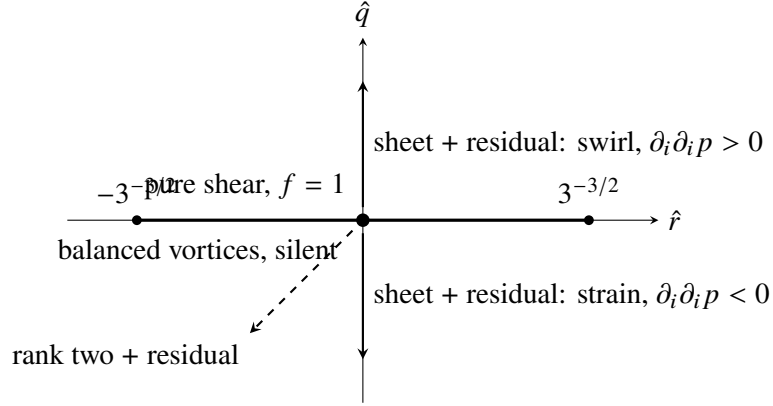

The first assumption decides the topology of the escape, not its source. For the generic rank-two nilpotent $\mathsf P$, the superposition of two layers, the cofactors do not vanish, since $\mathrm{adj}\,\mathsf P=\mathsf P^2$, so $R_a=-g^2\epsilon\,\tr(\mathsf P^2\boldsymbol b)+O(g\epsilon^2)$. For the two layers $u_1=x_2$ and $u_2=x_3$ the element is $b_{31}$, the residual's velocity along the second layer's normal varying along the first layer's velocity direction, which closes the chain of the two layers. The cubic is then dominated by its constant term \citep{moro_lidskii_1997}, and its three roots have magnitude $(g^2\epsilon)^{1/3}$, larger than $(g\epsilon)^{1/2}$ by $(g/\epsilon)^{1/6}$. But they sit at $120^\circ$, so their squares cancel. The escaping eigenvalues of the dressed rank-two nilpotent are balanced. Its source, $-2g\epsilon\,\tr(\mathsf P\boldsymbol b)$, is of the same order as the sheet's, and it is carried by the same gates, one for each layer, rather than by those eigenvalues. It leaves the origin with $|\hat q|$ and $|\hat r|$ both of order $\epsilon/g$, where the sheet has $|\hat r|$ of order $(\epsilon/g)^2$ (figure~\ref{fig:plane}). So the gate, the linearity of the source in the sheet strength and the sign of the source hold for either silent object. What depends on the rank is the eigenvalue geometry: the square-root or cube-root escape, the swirl-or-strain topology, and the active branch of $m$, which is $\hat q$ for the sheet and $\hat r$ for the pair of layers. Which branch is active on the near-cancelling set is therefore a single-point test of the first assumption.

\subsection{Leaving the sheet without leaving silence}\label{sec:ridge}

There is one escape that is not one. A gradient can leave the sheet and acquire eigenvalues of any size, yet stay silent, by growing the equilateral spectrum of \S\ref{sec:exception}. Along that ridge the swirl and the strain of the balanced vortex cancel exactly. Silence breaks only when the spectrum leaves the ridge, $Q_a\neq0$. Setting $Q_a=0$ in \eqref{eq:fissure} puts the ridge at $|n_ib_{ij}t_j|\le\epsilon/2g$. For a positive coupling the pair \eqref{eq:escape} sets the topology only where the second invariant dominates the cubic, $|n_ib_{ij}t_j|\gg(\epsilon/g)^{1/3}$. Between the two thresholds the sign of the source is already set by the coupling, but the discriminant of the cubic decides between a real and a complex pair. For a negative coupling there is no second threshold. For $g\gg\epsilon$ an open gate already leaves the ridge, so the ridge matters only where the residual is no longer small. The wall shuts the gate kinematically. There $\boldsymbol n=\boldsymbol e_2$ and $\boldsymbol t$ lies along $\boldsymbol\tau$, so the gate element is $\dd_1u_2$ taken along the wall stress, which no slip and continuity make $O(y^2)$.

\section{Discussion and consequences}\label{sec:conclusion}

The algebra of \S\S\ref{sec:bounds}--\ref{sec:escape} settles two things at every point of every incompressible flow. The source sees only the eigenvalues. Near-cancellation with a small third invariant is shear, to within the wedge's factor of three, and otherwise it is a balanced vortex. The sheet case settles a third. A dressed rank-one sheet escapes cancellation through one matrix element of its residual, at an amplitude set by the geometric mean of sheet and residual, with the sign of the coupling deciding between swirl and strain. At the wall the algebra is kinematic. No slip makes the wall plane a silent sheet with no escape at leading order, and its intensity is the wall shear-stress fluctuation. Whether turbulence uses this geometry is a hypothesis. It rests on two assumptions: that the silent object away from the wall is the rank-one sheet rather than the generic rank-two nilpotent, and that the residual is small beside it. The source rule does not rest on the first. The gate, the linearity in the sheet strength and the sign of the source hold for either silent object, and the first assumption fixes the topology of the escape.

The reference is a Gaussian random gradient with the same second moments, conditioned in the same way on intensity and near-cancellation, in the manner of the constrained random tensors of \citet{keylock_synthetic_2017}. Against it, the intense near-cancelling gradients of turbulence should be shear-dominated more often. They should be closer to rank one, with $r_1$ small, where $r_1/(1-r_1)=(\epsilon/g)^2$ when the sheet is the leading singular pair of $\aT$. Their vorticity should align with the intermediate strain axis more tightly than conditioning alone enforces. \citet{blackburn_topology_1996} found that alignment strengthening towards the wall without conditioning. The second-invariant branch of $m$ should be the active one, as it is for a dressed sheet, and the third invariant will say how much of the cancelling intensity balanced vortices carry instead.
Their geometry should be fixed in wall units above the buffer layer. \citet{arun_velocity_2024} report such a collapse for the mean partition of the gradient into normal and shearing parts, in boundary layers and to a lesser extent in a channel. The prediction here is for the conditioned set. Whether it is also the geometry of isotropic turbulence is a separate question, and needs a matched scaling. Each premise could fail.

Two consequences follow. The first is a proposal for the wall pressure. The escape algebra holds at every point at every Reynolds number. If the local geometry of the sheets is fixed in wall units, the growth of the inner-scaled wall-pressure variance \citep{klewicki_statistical_2008,panton_correlation_2017} cannot be a change in the rule by which the source escapes. It must be a change in how much escapes, and where. The sheets of the inertial layer are the vortical fissures. The velocity jump across a fissure is of the order of the friction velocity at every $\dt$ \citep{de_silva_interfaces_2017}, and the number of fissures across the layer grows as $\ln\dt$ \citep{de_silva_uniform_2016}. We have argued elsewhere that the residual a fissure carries scales in the same way, that the source of one fissure, integrated across its thickness, depends on the jump and the residual and not on the thickness, and that the fissures summed against the Poisson kernel give a logarithmic growth with a coefficient consistent with the measured one \citep{massey_behind_2026}. Equal delivery from every fissure needs one more thing. The kernel attenuates a source at height $y$ and wall-parallel wavenumber magnitude $k$ by $\mathrm e^{-ky}$. A higher fissure delivers as much as a lower one only if its along-sheet scale grows linearly with height, so that $ky$ stays of order one, while its local geometry stays fixed. The zones themselves thicken with distance from the wall \citep{de_silva_uniform_2016}. The along-sheet step is ours, and it is part of the proposal. The proposal is growth by extent: more sheets escaping at the same rate. Two Reynolds numbers cannot separate extent from amplitude, and leave open whether the sheets are also louder. Those tests check whether the local geometry is fixed in wall units. The delivery is a different question. The wall-pressure variance depends on the wall-parallel scale of each sheet's source and on whether the sources of different sheets add coherently, and neither is a property of one point. The net-source framework of \citet{anantharamu_analysis_2020} is the tool for that question.

The second consequence is for measurement and for models. The wedge is two-sided with constants, and exact at the wall. The shear fraction is therefore a pointwise proxy for cancellation by shear, and a robust one, since it does not require the small difference of two large terms. It says nothing about the balanced-vortex family, which the third invariant identifies. A model of the source, whether in a Reynolds-stress closure or for the subgrid pressure of a large-eddy simulation, must be built on what escapes. The source is the gate element times the sheet strength, whichever silent object carries it. Its topology depends on the rank. The largest gradients in the flow, the sheets themselves, contribute nothing to it however intense they are.

\section*{Funding}
The support of DARPA under award \#\,HR0011-24-9-0465 is gratefully acknowledged.

\section*{Declaration of interests}
The author reports no conflict of interest.

\section*{Data availability statement}
The script generating the random-tensor ensemble of appendix~\ref{app:tests} is available from the author on request.

\section*{Author ORCID}
J.M.O. Massey, \url{https://orcid.org/0000-0002-2893-955X}

\section*{AI statement}
Claude's Fable 5.1 model was used to check the algebra of \S\S2--4 and to refine the wording. The physical question, the interpretation of the results and the claims are the author's, who takes full responsibility.

\FloatBarrier
\appendix
\section{The wedge on random tensors}\label{app:tests}

\begin{figure}[htp]
    \centering
    \includegraphics[width=0.44\textwidth]{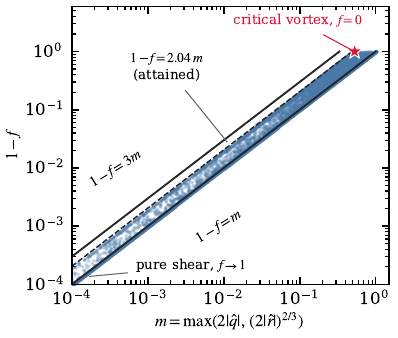}
    \caption{The wedge \eqref{eq:wedge} on $5.5\times10^5$ random trace-free tensors, Gaussian and near-shear. Every one lies between $1-f=m$ and $1-f=3m$, and none above the dashed line $1-f=2.04\,m$. Pure shear is at the origin. The critically strained vortex of \S\ref{sec:exception} is at $(m,1{-}f)=((2/3\sqrt3)^{2/3},1)$.}
    \label{fig:wedge}
\end{figure}

Figure~\ref{fig:wedge} shows the wedge \eqref{eq:wedge} on $5.5\times10^5$ random trace-free tensors. The ensemble has two parts. In the larger, $4\times10^5$ tensors, every component is drawn independently from a Gaussian and the trace is removed. In the smaller, $1.5\times10^5$ tensors, a unit rank-one layer $\boldsymbol t\otimes\boldsymbol n$ is dressed with a random trace-free residual of norm $\epsilon$, with $\epsilon$ spread over four decades, which fills the neighbourhood of the origin that the Gaussian part leaves empty. Every point lies between $1-f=m$ and $1-f=3m$, and none above $1-f=2.04\,m$. The balanced-vortex family, with the critically strained vortex at its end, lies on the line $1-f=3\cdot2^{-2/3}m$, inside the wedge, since on it $1-f=3|\hat r|^{2/3}$ and $m=(2|\hat r|)^{2/3}$. Every identity and bound of \S\S\ref{sec:identity}--\ref{sec:escape} was checked on the same ensemble.

\FloatBarrier

\bibliographystyle{plainnat}
\bibliography{refs}

\end{document}